\documentclass[12pt]{article}
\usepackage{amsmath,amsthm,amsfonts,amssymb,amscd,amsbsy}
\usepackage[latin1]{inputenc}

\newcommand{\la}{{\lambda}}

\newcommand{\re}{{\mathbb{R}}}

\newcommand{\N}{{\mathbb{N}}}

\newcommand{\al}{{\alpha}}

\newcommand{\ve}{{\varepsilon}}

\def\ov{\overline}

\begin{document}

\centerline{\Large\bf There are no $C^1$-stable intersections of}

\vskip .1in

\centerline{\Large\bf regular Cantor sets}

\vskip .2in

\centerline{\sc Carlos Gustavo Moreira}

\centerline{\sl Instituto de Matem\'atica Pura e Aplicada}

\centerline{\sl Estrada Dona Castorina 110}

\centerline{\sl 22460-320 Rio de Janeiro, RJ}

\centerline{\sl Brasil}

\vskip .4in

\noindent
{\bf 1. Introduction:} The existence of stable intersections of regular Cantor sets is a fundamental tool to provide persistent examples of non-hyperbolic $C^2$ diffeomorphisms of surfaces, as did Newhouse ([N]), by means of the concept of thickness of a Cantor set. The thickness is a fractal invariant, which is continuous and positive in the $C^2$ (or even in the $C^{1+\al}$, $0<\al<1)$ topology, such that, if the product of the thicknesses of two regular Cantor sets is larger than one and their support intervals intersect is a nontrivial way, then they have stable intersection (see [PT]). However Ures ([U]) showed that, in the $C^1$ topology, the  thickness of regular Cantor sets is terribly discontinuous. Indeed, generic $C^1$-regular Cantor sets have zero thickness. He also showed that two regular Cantor sets whose support intervals touch at one point cannot have extremal stable intersection (in the sense of [M]) in the $C^1$ topology and he used these results to show that $C^1$-generic first homoclinic bifurcations present full upper density of hyperbolicity at the initial bifurcation parameter.

However, despite the discontinuity of the thickness in the $C^1$ topology, the Hausdorff dimension of regular Cantor sets is continuous and positive in the $C^1$ topology (and coincides with the limit capacity). On the other hand, it was showed in [MY] that generic pairs of regular Cantor sets in the $C^2$ (or $C^{1+\al}$) topology whose sum of Hausdorff dimensions is larger than one have translations which have stable intersection. Moreover, they have translations whose intersections have stably positive Hausdorff dimensions. This poses a more difficult problem: is it always possible to destroy intersections of regular Cantor sets by performing arbitrarily small $C^1$ perturbations of them? The situation is particularly delicate when the intersection between the Cantor sets has positive Hausdorff dimension, which is a typical situation in the $C^2$ topology as seen before. We solve this problem in the following

\vskip .2in

\noindent
{\bf Theorem 1.} {\it Given any pair $(K,K')$ of  regular Cantor sets, we can find, arbitrarily close to it in the $C^1$ topology, pairs $(\tilde K,\tilde K')$ of regular Cantor sets with $\tilde K\cap \tilde K'=\emptyset$.

Moreover, for generic pairs $(K,K')$ of $C^1$-regular Cantor sets, the arithmetic difference $K-K'=\{x-y \mid x\in K, y\in K'\} = \{t\in\re \mid K\cap(K'+t)\ne\emptyset\}$ has empty interior (and so is a Cantor set).}

This answers a question by Christian Bonatti.

Since stable intersections of Cantor sets are the main known obstructions to density of hyperbolicity for diffeomorphisms of surfaces, the previous result gives some hope of proving density of hyperbolicity in the $C^1$ topology for diffeomorphisms of surfaces. In particular it is used in a forthcoming joint work with Carlos Matheus and Enrique Pujals on a family of two-dimensional maps (the so-called Benedicks-Carleson {\it toy model} for H\'enon dynamics) which present stable homoclinic tangencies (Newhouse's phenomenon) in the $C^2$-topology but whose elements can be arbitrarily well approximated in the $C^1$-topology by hyperbolic maps. 

The main technical difference between the $C^1$ case and the $C^2$ (or even $C^{1+\al}$) cases is the lack of bounded distortion of the iterates of $\psi$ in the $C^1$ case, and this fact will be fundamental for the proof of the previous result.

I would like to thank Carlos Matheus for helpful comments and suggestions which substantially improved this work.

\vskip .2in

\noindent
{\bf 2. Proofs of the results:}

We recall that $K$ is a {\it $C^k$-regular Cantor set}, $k\ge 1$, if:

\begin{itemize}

\item[i)] there are disjoint compact intervals $I_1,I_2,\dots,I_r$ such that $K\subset I_1 \cup \cdots\cup I_r$ and the boundary of each $I_j$ is contained in $K$;

\item[ii)] there is a $C^k$ expanding map $\psi$ defined in a neighbourhood of $I_1\cup I_2\cup\cdots\cup I_r$ such that $\psi(I_j)$ is the convex hull of a finite union of some intervals $I_s$ satisfying:

\begin{itemize}

\item[ii.1)] for each $j$, $1\le j\le r$ and $n$ sufficiently big, $\psi^n(K\cap I_j)=K$;

\item[ii.2)] $K=\bigcap\limits_{n\in\mathbb N} \psi^{-n}(I_1\cup I_2\cup\cdots\cup I_r)$.

\end{itemize}
\end{itemize}

We say that $\{I_1,I_2,\dots,I_r\}$ is a {\it Markov partition} for $K$ and that $K$ is {\it defined} by $\psi$.

Given $s \in [1,k]$ and another regular Cantor set $\tilde K$, we say that $\tilde K$ is close to $K$ in the $C^s$ topology if $\tilde K$ has a Markov partition $\{\tilde I_1,\tilde I_2,\dots,\tilde I_r\}$ such that the interval $\tilde I_j$ has endpoints close to the endpoints of $I_j$, for $1 \le j \le r$ and $\tilde K$ is defined by a $C^s$ map $\tilde \psi$ which is close to $\psi$ in the $C^s$ topology.

Given a $C^{3/2}$-regular Cantor set $K$, we define the parameter $\la(K)=\max\{|\psi'(x)|, x\in \bigcup\limits_{j=1}^r I_j\}>0$, which depends continuously on $K$ in the $C^1$-topology.

We may associate for each $j\le r$ a gap $U_j\subset I_j$ of $K$, which is determined by the combinatorics of $(K, \psi)$ in the following way: we take the smallest $m_j \ge 1$ such that $\psi^{m_j}(I_j)$ is the convex hull of the union of more than one interval of the Markov partition of $K$, say of $I_{s_j} \cup I_{s_j+1} \cup \dots \cup I_{s_j+l_j}$, and, if $V_j$ is the gap between $I_{s_j}$ and $I_{s_j+1}$, we set $U_j:=\psi^{-m_j}(V_j)$. We define a parameter $b(K)$ in the following way: given $n\in\mathbb N$, $j\le r$ and $\tilde J$ a connected component of $\psi^{-n}(I_j)$, we define $b(\tilde J)=|U|/|\tilde J|$, where $U\subset\tilde J$ is the gap of $K$ such that $\psi^n(U)=U_j$, and we define $b(K)=\underset{\tilde J}{\text{inf}}\{b(\tilde J)$; $\tilde J$ connected component of $\psi^{-n}(I_j)$, $\exists\, n\in\mathbb N$, $j\le r \}>0$. We also define a similar parameter $g(K)$ as follows: given $n\in\mathbb N$, $j\le r$ and $\tilde J$ a connected component of $\psi^{-n}(I_j)$, we define $L_{\tilde J}$ and $R_{\tilde J}$ the gaps of $K$ attached to the left and right endpoints of $\tilde J$, respectively. We set $g(\tilde J)=\min\{|L_{\tilde J}|/|\tilde J|,|R_{\tilde J}|/|\tilde J|\}$, and we define $g(K)=\underset{\tilde J}{\text{inf}}\{g(\tilde J)$; $\tilde J$ connected component of $\psi^{-n}(I_j)$, $\exists\, n\in\mathbb N$, $j\le r \}>0$. Finally, we define $a(K)=\min\{b(K),g(K)\}$. Notice that $b(K), g(K)$ and $a(K)$ depend continuously on $K$ in the $C^{3/2}$ topology.

In the next lemma we will exploit the lack of bounded distortion in order to produce a very distorted geometry (with very large gaps) near some subsets of a $C^1$-regular Cantor set.

\vskip .2in

\noindent
{\bf Lemma 1.} {\it Let $K$ be a $C^2$-regular Cantor set. Let $c(K)=2\la(K)/a(K)$. Then, given $\ve>0$, let $n_0=\lceil c(K)\log \ve^{-1}/\ve \rceil$. Suppose that $X\subset K$ is a compact set satisfying $\psi^i(X)\cap \psi^j(X)=\emptyset$, for $0\le i<j\le n_0$. Then, for any $\delta>0$,
we can find a covering of $K$ formed by intervals $(J_i)$ of its construction (i.e. intervals which are connected components of $\psi^{-n}(I_j)$ for some $n\in\mathbb N$, $j\le r$) which have size smaller than $\delta$ satisfying the following properties: 

Let $D$ be the union of the intervals $J_i$ for all $i$ and the intervals $\psi^j(J_i), 1 \le j \le n_0$, for the intervals $J_i$ which intersect $X$. There is a Cantor set $\tilde K$ in the $\ve$-neighbourhood of $K$ in the $C^1$-topology with $a(\tilde K) \ge a(K)$ such that all connected components of $D$ are still intervals of the construction of $\tilde K$ and $J_i\cap X \ne \emptyset \Rightarrow \tilde K\cap J_i$ has a gap $V_i$ with $|V_i|\ge (1-\ve)|J_i|$\/}.

\vskip .2in

\noindent
{\bf Proof:} Let $A=\{i\mid J_i\cap X\ne\emptyset\}$, where $\{J_i\}$ is the convering of $K$ by the connected components of $\psi^{-N}(\bigcup\limits_{j=1}^r I_j)$, for some large $N$ so that, in particular, $|J_i|<\delta$, for all $i$.

Since $X$ is compact, if $N$ is large enough, we have by the hypothesis of the lemma, $\psi^i (\tilde X)\cap \psi^j(\tilde X)=\emptyset$ for $0\le i<j\le n_0$, where $\tilde X:=\bigcup\limits_{i\in A} J_i$. 

We may perform, as in Lemma II.2.1 of [M], a small change in $\psi$ in the $C^{3/2}$ topology in such a way that the restrictions of $\psi$ to the intervals $\psi^j(J_i)$ with $i\in A$, $0\le j< n_0$ become affine; we change $\psi$ just in these intervals and in the gaps attached to them.

Now we will make small $C^1$ perturbations on the restriction of $\psi$ to the intervals $\psi^j(J_i)$ with $i\in A$, $0\le j<n_0$. We will begin changing $\psi$ in the intervals $\psi^{n_0-1}(J_i)$, then in the intervals $\psi^{n_0-2}(J_i)$ and so on, in order to make the proposition of a gap in each of these intervals grow in such a may that the size of each of the two remaining intervals is multiplied by $1-\frac{2a\ve}{3\la}$, where $\la=\la(K)$ and $a=a(K)$.

More precisely, if $\psi^j(J_i)=[r,s]$, $i\in A$, $0\le j< n_0$ is some of these intervals and $m$ is such that $\psi^m(J_i)=I_l$, let $\tilde U\subset J_i$ such that $\psi^m(\tilde U)=U_l$ and let $\psi^j(\tilde U)=(u,v)\subset[r,s]$. Writing $\psi|_{[r,s]}(x)=\tilde \la x + t$, we consider the affine map $\tilde \psi|_{[r,s]}$ given by 
$$
\tilde\psi|_{[r,s]}(x) = \begin{cases}
\tilde\la (1-\frac{2 a\ve}{3\la})^{-1} (x-r)+\tilde \la r+t,\text{ if } x\in[r,r+(1-\frac{2a\ve}{3\la})(u-r)] \\
\tilde\la (1-\frac{2 a\ve}{3\la})^{-1} (x-s)+\tilde \la s+t,\text{ if } x\in[s-(1-\frac{2a\ve}{3\la})(s-v),s]
\end{cases}
$$
and we extend $\tilde \psi$ to $[r,s]$ in such a way that $\tilde\psi|_{[r,s]}$ is a $C^1$ function. Notice that the image $\tilde\psi((r+(1-\frac{2a\ve}{3\la})(u-r),s-(1-\frac{2a\ve}{3\la})(s-v)))=\psi((u,v))$ of the gap remains the same. The size of the new gap in $[r,s]$ is $v-u+\frac{2a\ve}{3\la}(s-v+u-r)<v-u+\frac{2a\ve}{3\la}(s-r)<(1+\frac{2\ve}{3\la})(v-u)$. In particular, it is not difficult to see that we may construct such a function $\tilde\psi$ with $||\tilde\psi-\psi||_{C^1}<\ve$.

Finally, the total proportion of the complement of the new gap $V_i=\tilde U$ for the modified $\psi$ (indeed $\tilde\psi|_{J_i}^{-n_0}(\psi^{n_0}(\tilde U)))$ is at most
$(1-\frac{2a\ve}{3\la})^{n_0} (1-a) \le(1-\frac{2a\ve}{3\la})^{\frac{2\la}{a\ve}\log \ve^{-1}} \cdot (1-a) < \ve^{4/3}<\ve$. It is not difficult to see that after these perturbations we will also have $a(\tilde K) \ge a(K)$ (indeed in the non-affine part of the dynamics we are only increasing  the proportion of some gaps, and in the affine and local part of the dynamics the proportion of the gaps is preserved), which concludes the proof of the Lemma. $\qed$

\vskip .3in

\noindent
{\bf Lemma 2.} {\it Given a $C^2$-regular Cantor set $K'$, for a residual set of  $C^2$-regular Cantor sets $K$, if $k= \lfloor{(1- HD(K'))}^{-1}\rfloor +1$ then $\bigcap \limits_{j=1}^{k}A_j = \emptyset$ , where $A_j=F_j((K \cap K') \cap P_j), 1 \le j \le k$ and $(F_j, P_j), 1 \le j \le k$ are distinct elements of \{$(\psi^r|_I,I); r \in \N$ and $I$ is a maximal interval of the construction of $K$ where $\psi^r$ is injective\}.}

\vskip .2in

\noindent
{\bf Proof:}  It is enough to show that, for generic $C^2$-regular Cantor sets $K$, given distinct sets $A_j=\psi^{r_j}((K \cap K') \cap P_j), 1 \le j \le k$, where $P_j$ is a maximal interval of the construction of $K$ such that $\psi^{r_j}|_{P_j}$ is injective, we have $\bigcap \limits_{j=1}^{k}\psi^{r_j}(K \cap K' \cap P_j) = \emptyset$. So, we will, from now on, fix the branches $\psi^{r_j}|_{P_j}$ that we will consider.

Typically, $K \cap K'$ does not contain any preperiodic point of $\psi$ (indeed, there is only a countable number of them in $K$, and so almost all translations of $K'$ do not intersect them). In this case, given a point $x \in \bigcap \limits_{j=1}^{k}A_j$, with $x=\psi^{r_j}(y_j), 1 \le j \le k$ with $y_j \in K \cap K' \cap P_j, \forall j \le k$, we have that the points $y_j$ are all distinct (indeed, if $y_i=y_j$ with $i \neq j$ then $r_i \neq r_j$, since $A_i$ and $A_j$ are distinct, so $x$ is periodic and the points $y_i, y_j \in K \cap K'$ are preperiodic). We may find disjoint intervals $\hat{P_j}(x) \subset P_j$ of the construction of $K$ with $y_j \in \hat{P_j}(x), 1 \le j \le k$, and take $\ve_x>0$ and a neighbourhood ${\cal V}_x$ of $K$ in the $C^1$ topology such that, for $\tilde K \in {\cal V}_x$, if $\tilde P_j$, $\tilde {\hat {P_j}}(x)$ and $\tilde \psi$ denote the continuations of the intervals $P_j, \hat{P_j}(x)$ and of the map $\psi$ which defines $K$ respectively, we have $\bigcap \limits_{j=1}^{k} \tilde \psi^{r_j}((\tilde K \cap K')\cap \tilde P_j) \cap (x-\ve_x, x+\ve_x) \subset \bigcap \limits_{j=1}^{k} \tilde \psi^{r_j}((\tilde K \cap K')\cap \tilde {\hat{P_j}}(x))$.  We may take a finite covering of the compact set $\bigcap \limits_{j=1}^{k}A_j$ by intervals $(x_i-\ve_{x_i}, x_i+\ve_{x_i}), 1 \le i \le m$. We denote by $\tilde {\cal V}$ the neighbourhood $\bigcap \limits_{i=1}^{m} {\cal V}_{x_i}$ of $K$ in the $C^1$ topology. We will assume (by reducing $\tilde {\cal V}$, if necessary) that, for any $\tilde K \in \tilde {\cal V}$,  $\bigcap \limits_{j=1}^{k} \tilde \psi^{r_j}((\tilde K \cap K')\cap \tilde P_j) \subset \bigcup \limits_{i=1}^{m} (x_i-\ve_{x_i}, x_i+\ve_{x_i})$, so $\bigcap \limits_{j=1}^{k} \tilde \psi^{r_j}((\tilde K \cap K')\cap \tilde P_j) \subset \bigcup \limits_{i=1}^{m} \bigcap \limits_{j=1}^{k} \tilde \psi^{r_j}((\tilde K \cap K')\cap \tilde {\hat{P_j}}(x_i))$.  It is enough to show that, for each $i \le m$, there is an open and dense set $\tilde {\cal V}_i \subset \tilde {\cal V}$ such that, for any $\tilde K \in \tilde {\cal V}_i$, we have $\bigcap \limits_{j=1}^{k} \tilde \psi^{r_j}((\tilde K \cap K')\cap \tilde {\hat{P_j}}(x_i))=\emptyset$. Since the sets $\psi^{r_j}((\tilde K \cap K')\cap \tilde {\hat{P_j}}(x_i))$ are compact, the above condition is clearly open, so it is enough to prove that it is dense in $\tilde {\cal V}$.

Now, since the intervals $\tilde {\hat{P_j}}(x_i), 1 \le j \le k$ are disjoint, we may consider families of perturbations ${\tilde K}_{t_1,\cdots,t_k} \in  \tilde {\cal V}$ of $K$ defined by maps ${\tilde \psi}_{t_1,\cdots,t_k}$ which form a family of perturbations of $\tilde \psi$ depending on $k$ small parameters $t_1, t_2, \cdots, t_k \in (-\delta, \delta)$, for some (small) $\delta>0$ for which ${\tilde \psi^{r_j}}_{t_1,\cdots,t_k}|_{\tilde{\hat{P_j}}(x_i)}= {\tilde \psi^{r_j}}|_{\tilde{\hat{P_j}}} + t_j, \forall t_1, t_2, \cdots, t_k \in (-\delta, \delta)$ (it is enough to make suitable perturbations of $\tilde \psi$ in the disjoint intervals $\tilde {\hat{P_j}}(x_i), 1 \le j \le k$). Since the limit capacities (box dimensions) of the sets ${\tilde \psi}^{r_j}(K'\cap \tilde P_j)$ are bounded by $HD(K')$, $\prod \limits_{j=1}^k \tilde \psi^{r_j}(K'\cap \tilde P_j)$ has limit capacity bounded by $k.HD(K')$$\ < k - 1$, and so linear projections of it in ${\mathbb R}^{k - 1}$ have zero Lebesgue measure.
Since \{$(t_1,\cdots,t_k)\in {\mathbb R}^k \mid \bigcap \limits_{j=1}^k (\psi^{r_j}(K' \cap \tilde P_j) + t_j)\neq \emptyset$\} = $\{(t, t+ s_1, \cdots, t+ s_{k-1}) \mid (s_1, \cdots, s_{k- 1})= D_k(\prod \limits_{j= 1}^k \psi^{r_j}(K' \cap \tilde P_j)) \}$, where $D_k$ is the linear map given by $D_k(x_1, \cdots, x_k):= (x_2 - x_1, x_3 - x_1, \cdots, x_k - x_1)$, for almost all $(t_1, \cdots, t_k)$ the intersection is empty, which implies the result.   $\qed$

\vskip .2in
\noindent
{\bf Lemma 3.} {\it Let $(K,K')$ be a pair of $C^2$-regular Cantor sets and $B= \{(\psi^r|_I,I); r \in \N$ and $I$ is a maximal interval of the construction of $K$ where $\psi^r$ is injective\}. Given any $m \ge 1$, the following holds: for any fixed distinct elements $(\psi^{r_j}|_{P_j},P_j), 1 \le j \le m$ of $B$ and for any $\eta>0$ there is a $C^2$-regular Cantor set $\tilde K$ at a distance smaller than $\eta$ from $K$ in the $C^1$ topology and with $a(\tilde K)>a(K)-\eta$ such that $\bigcap \limits_{j=1}^m {\tilde \psi}^{r_j}(\tilde K \cap K' \cap P_j)= \emptyset$, where $\tilde \psi$ is the map which defines $\tilde K$.} 

\vskip .2in
\noindent
{\bf Proof:} Since the quantity $a(K)$ depends continuously on $K$ in the $C^{3/2}$-topology (and, in particular, in the $C^2$-topology), Lemma 2 implies that the above statement is true for $m \ge k= \lfloor{(1- HD(K'))}^{-1}\rfloor +1$. We will argue by backward induction on $m$: we will show that, if $q \ge 1$ and the statement of Lemma 3 is true for every $m \ge q+1$ then it is true for $q$.

Let $(K, K')$ be a pair of $C^2$-regular Cantor sets, $\eta  \in (0,1)$ and $(\psi^{r_j}|_{P_j},P_j), 1 \le j \le q$ fixed distinct elements of $B$. Let $X = \bigcap \limits_{j=1}^q \psi^{r_j}(K \cap K' \cap P_j)$. Since $K'$ is a $C^2$-regular Cantor set, it has bounded geometry; in particular, there is $\ov \lambda >0$ such that, for every interval $J'_1$ of the construction of $K'$, there is an interval $J'_0$ of a previous step of the construction of $K'$ which strictly contains $J'_1$ such that $|J'_0|/|J'_1|<\ov \lambda$. Let $\ve=a(K) a(K') \eta/5 \ov \la$ and let $N_0=\lceil 2 c(K)\log \ve^{-1}/\ve \rceil$, where $c(K)=2\la(K)/a(K)$. For each pair $(i,j)$ with $0 \le i < j \le N_0$, the intersection $\psi^i(X) \cap \psi^j(X)$ can be written as a finite union of intersections of at least $q+1$ sets of the form $\psi^r(K \cap K' \cap I)$, where $(\psi^r|_I,I) \in B$. By the induction hypothesis (applied several times to make a sequence of small perturbations of the first Cantor set, one for which one of the intersections mentioned above), we may approximate $K$ by a $C^2$-regular Cantor set $\check K$ at a distance smaller than $\eta/2$ from $K$ in the $C^1$ topology with $a(\check K)>\max\{a(K)-\eta/2, a(K)/2\}$ and $c(\check K)<4\la(K)/a(K)$,  such that, if $\check X = \bigcap \limits_{j=1}^q  {\check \psi}^{r_j}(\check K \cap K' \cap {\check P_j})$, where $\check \psi$ is the map which defines $\check K$, then the sets ${\check \psi}^j(\check X), 0 \le j \le N_0$ are pairwise disjoint. So, if $\check Y:=({\check \psi}^{r_1}|_{\check P_1})^{-1}(\check X)=(\check K \cap K' \cap {\check P_1})\cap\bigcap \limits_{j=2}^q  ({\check \psi}^{r_1}|_{\check P_1})^{-1}({\check \psi}^{r_j}(\check K \cap K' \cap {\check P_j}))$, then the sets ${\check \psi}^j(\check Y), 0 \le j \le N_0$ are pairwise disjoint. Indeed, if ${\check \psi}^i(\check Y) \cap {\check \psi}^j(\check Y) \neq \emptyset$, with $0 \le i < j \le N_0$, we would have $\emptyset \neq {\check \psi}^{r_1}({\check \psi}^i(\check Y)) \cap {\check \psi}^{r_1}({\check \psi}^j(\check Y))={\check \psi}^i({\check \psi}^{r_1}(\check Y)) \cap {\check \psi}^j({\check \psi}^{r_1}(\check Y))={\check \psi}^i(\check X) \cap {\check \psi}^j(\check X)=\emptyset$, a contradiction. This implies, by compacity, that we may find a covering of $\check Y$ by (small) intervals $\check J^*_i$ of the construction of $\check K$ such that, if $Y^* \supset \check Y$ is the intersection of $\check K$ with the union of the intervals  $\check J^*_i$, we still have that the sets ${\check \psi}^j(Y^*), 0 \le j \le N_0$ are pairwise disjoint.

So, $Y^*$ satisfies the hypothesis of Lemma 1 for $\ve$ (since $N_0=\lceil 2 c(K)\log \ve^{-1}/\ve \rceil\ge \break \lceil c(\check K)\log \ve^{-1}/\ve \rceil$), and thus the conclusion of Lemma 1 holds: 
for any $\delta>0$, we can find a covering of $\check K$ formed by intervals $(\check J_i)$ of its construction which have size smaller than $\delta$ satisfying the following properties: let $D$ be the union of the intervals $\check J_i$ for all $i$ and the intervals ${\check \psi}^j(\check J_i), 1 \le j \le n_0$, for the intervals $\check J_i$ which intersect $Y^*$. There is a Cantor set $\ov K$ in the $\ve$-neighbourhood of $\check K$ in the $C^1$-topology with $a(\ov K) \ge a(\check K)$ such that all connected components of $D$ are still intervals of the construction of $\ov K$ and $\check J_i\cap Y^* \ne \emptyset \Rightarrow \ov K\cap \check J_i$ has a gap $V_i$ with $|V_i|\ge (1-\ve)|\check J_i|$\/. If $\delta$ is small enough, the intervals $\check J^*_i$ are still intervals of the construction of $\ov K$, and $(\ov K \cap K' \cap {\ov P_1})\cap\bigcap \limits_{j=2}^q  ({\ov \psi}^{r_1}|_{\ov P_1})^{-1}({\ov \psi}^{r_j}(\ov K \cap K' \cap {\ov P_j}))$ is contained in the union of the intervals  $\check J^*_i$, and so is contained in the union of the intervals $\check J_i$ which intersect $Y^*$.

Now, we can apply a $C^1$ diffeomorphism $(\eta/2)$-close to the identity in the $C^1$ topology to $\ov K$ in order to make it disjoint from $K'$. Indeed, let $(\check J_i)^{(1)}$ and $(\check J_i)^{(2)}$ be the connected components of $\check J_i \setminus V_i$. We will make small independent translations of these intervals (if they do intersect $K'$) in the following way: if such an interval $(\check J_i)^{(s)}$ intersects $K'$, take an interval $J'$ of the construction of $K'$ intersecting it whose size belongs to the interval $(|(\check J_i)^{(s)}|/a(K'), \ov \la |(\check J_i)^{(s)}|/a(K')]$; the gaps attached to the ends of this interval have size larger than $|(\check J_i)^{(s)}|$ so we can apply a translation of it of size at most $\ov \la |(\check J_i)^{(s)}|/a(K')$ whose image is contained in one of these gaps. These translations can be performed all together (for all $i, s$) by a diffeomorphism at a $C^1$ distance to the identity smaller than $\eta/2$, since the gaps attached to the intervals $(\check J_i)^{(s)}$ have size at least $a(\check K)|J_i|$ and the size of the translations is at most $\ov \la |(\check J_i)^{(s)}|/a(K')<\ve \ov \la |\check {J_i}|/a(K')=a(\check K)|\check {J_i}|\eta/5$. If we denote by $\tilde K$ the image of $\check K$ by this diffeomorphism (which is a regular Cantor set defined by a map $\tilde \psi$ conjugated to $\ov \psi$ by the diffeomorphism that we applied to $\ov K$) we will have $(\tilde K \cap \check J_i) \cap K'$ empty, $\forall i$, so $\tilde Y:=(\tilde K \cap K' \cap {\tilde P_1})\cap\bigcap \limits_{j=2}^q  ({\tilde \psi}^{r_1}|_{\tilde P_1})^{-1}{\tilde \psi}^{r_j}(\tilde K \cap K' \cap {\tilde P_j}) =\emptyset$, and, applying ${\tilde \psi}^{r_1}$, we get \break $\tilde X := \bigcap \limits_{j=1}^q  {\tilde \psi}^{r_j}(\tilde K \cap K' \cap {\tilde P_j})={\tilde \psi}^{r_1}(\tilde Y)=\emptyset$, and we are done. $\qed$

\vskip .3in
\noindent
{\bf Proof of Theorem 1.} Lemma 3 for $m=1$ implies that generically $K \cap K'= \emptyset$. It follows, as in Theorem I.1 of [M], that, for each $r \in \mathbb Q$, \{$(K,K') \mid r \notin K - K'$\}= \{$(K,K') \mid K \cap (K'+ r) = \emptyset$\} is residual, and thus
\{$(K,K') \mid (K - K')\cap \mathbb Q = \emptyset$\} is residual. So, generically, for $(K, K')$ pair of $C^1$-regular Cantor sets, $K - K'$ has empty interior, and so is a Cantor set. $\qed$

\end{document}